\documentclass[12pt,leqno]{article}

\usepackage{amsmath,amssymb,amsthm}
\usepackage[all]{xy}
\usepackage{lscape}

\makeatletter

\newtheorem{theorem}{Theorem}

\newtheorem{lemma}[theorem]{Lemma}
\newtheorem{corollary}[theorem]{Corollary}

\theoremstyle{definition}

\theoremstyle{remark}

\theoremstyle{remark}

\def\({{\rm (}}
\def\){{\rm )}}

\let\Mathrm\operator@font
\let\Cal\mathcal
\let\Bbb\mathbb

\def\standop#1{\mathop{\Mathrm #1}\nolimits}
\def\difstop#1#2{\expandafter\def\csname #1\endcsname{\standop{#2}}}
\def\defstop#1{\difstop{#1}{#1}}

\defstop{AB}
\defstop{ann}
\defstop{Ass}

\defstop{CMFI}
\defstop{codim}
\defstop{Coh}
\defstop{Coker}
\defstop{Cone}
\defstop{Cl}

\defstop{depth}

\defstop{EM}
\defstop{embdim}
\defstop{End}
\defstop{ev}
\defstop{Ext}

\defstop{Flat}
\defstop{Func}

\difstop{height}{ht}
\defstop{Hom}

\def\id{\mathord{\Mathrm{id}}}

\difstop{Image}{Im}
\defstop{ind}

\defstop{Ker}

\defstop{Lch}
\defstop{length}
\defstop{Lin}
\defstop{Lqc}
\defstop{lqc}

\defstop{Mat}
\defstop{Max}
\defstop{Min}
\defstop{Mod}
\defstop{Mor}
\defstop{MCM}

\defstop{Nerve}
\defstop{NonSFR}
\defstop{NonCMFI}
\defstop{NonNor}
\defstop{Nor}

\defstop{PA}
\defstop{PM}
\defstop{Proj}

\defstop{Qch}
\defstop{qch}

\defstop{rank}

\defstop{res}
\defstop{Reg}

\defstop{Spec}
\defstop{supp}
\defstop{Supp}
\defstop{Sym}
\defstop{Sing}
\defstop{SFR}
\defstop{Soc}
\def\St{St}

\difstop{tdeg}{trans.deg}

\defstop{Tor}
\difstop{trace}{tr}

\defstop{Zar}

\def\O{\Cal O}

\def\sdarrow#1{\downarrow\hbox to 0pt{\scriptsize$#1$\hss}}
\def\suarrow#1{\uparrow\hbox to 0pt{\scriptsize$#1$\hss}}
\def\ssearrow#1{\searrow\hbox to 0pt{\scriptsize$#1$\hss}}

\def\section{\@startsection{section}{1}{\z@ }%
{-3.5ex plus -1ex minus -.2ex}{2.3ex plus .2ex}{\bf }}

\long\def\refname{\par\kern -3ex
\begin{center}\rm R\sc{eferences}\end{center}\par\kern 
-2ex}

\def\@seccntformat#1{\csname the#1\endcsname.\quad}

\def\@@@sect#1#2#3#4#5#6[#7]#8{%
   \ifnum #2>\c@secnumdepth 
      \def \@svsec {}\else \refstepcounter {#1}%
      \def\@svsec{}
   \fi 
   \@tempskipa #5\relax 
   \ifdim \@tempskipa >\z@ 
     \begingroup #6\relax \@hangfrom {\hskip #3\relax 
     \@svsec}{\interlinepenalty \@M #8\par }\endgroup 
     \csname #1mark\endcsname {#7}
   \else 
   \def \@svsechd {#6\hskip #3\@svsec #8\csname #1mark\endcsname {#7}}
   \fi \@xsect {#5}}

\def\@@@startsection#1#2#3#4#5#6{%
 \if@noskipsec \leavevmode \fi \par \@tempskipa #4\relax \@afterindenttrue 
 \ifdim \@tempskipa <\z@ \@tempskipa -\@tempskipa \@afterindentfalse 
 \fi \if@nobreak \everypar {}\else \addpenalty {\@secpenalty }\addvspace 
  {\@tempskipa }\fi \@ifstar {\@ssect {#3}{#4}{#5}{#6}}{\@dblarg 
  {\@@@sect {#1}{#2}{#3}{#4}{#5}{#6}}}}

\def\theparagraph{\thesection.\arabic{paragraph}}
\def\aparagraph{\@@@startsection{paragraph}{2}{\z@ }%
              {1.75ex plus .2ex minus .15ex}{-1em}{\bf(\theparagraph) } }
\def\paragraph{\@@@startsection{paragraph}{2}{\z@ }%
              {1.75ex plus .2ex minus .15ex}{-1em}{}{\bf(\theparagraph)} }

\let\c@theorem\c@paragraph

\title{Good filtrations and $F$-purity of invariant subrings}
\author{M{\sc itsuyasu} H{\sc ashimoto}}
\date{\normalsize
Graduate School of Mathematics, Nagoya University\\
Chikusa-ku,  Nagoya 464--8602 JAPAN\\
{\small \tt hasimoto@math.nagoya-u.ac.jp}}

\begin{document}

\maketitle

\begin{abstract}
Let $k$ be an algebraically closed field of positive characteristic, 
$G$ a reductive group over $k$, and $V$ a finite dimensional $G$-module.
Let $B$ be a Borel subgroup of $G$, and $U$ its unipotent radical.
We prove that if $S=\Sym V$ has a good filtration, then $S^U$ is $F$-pure.
\end{abstract}

Throughout this paper, $p$ denotes a prime number.
Let $k$ be an algebraically closed field of characteristic $p$, 
and $G$ a reductive group over $k$.
Let $B$ be a Borel subgroup of $G$, and $U$ its unipotent radical.
We fix a maximal torus $T$ contained in $B$, 
and fix a base of the root system $\Sigma$ of $G$ so that $B$ is negative.
For any weight $\lambda\in X(T)$, we denote the induced module
$\ind_B^G(\lambda)$ by $\nabla_G(\lambda)$.
We denote the set of dominant weights by $X^+$.
For $\lambda\in X^+$, we call $\nabla_G(\lambda)$ the dual Weyl module
of highest weight $\lambda$.
We say that a $G$-module $W$ has a good filtration \cite{Donkin}
if $H^1(G,W\otimes \nabla_G(\lambda))=0$ for any $\lambda\in X^+$.

Let $V$ be a finite dimensional $G$-module, and $S=\Sym V$.
The objective of this paper is to prove the following.

\begin{theorem}\label{main.thm}
If $S$ has a good filtration, then $S^U$ is $F$-pure.
\end{theorem}

For a commutative ring $R$ and an $R$-linear map $\varphi:M\rightarrow N$ of
$R$-modules, we say that $\varphi$ is pure or $R$-pure if
$1_W\otimes \varphi:W\otimes_R M\rightarrow W\otimes_R N$ is injective
for any $R$-module $W$.
We say that $\varphi$ is a split mono if there is an $R$-linear map
$\psi:N\rightarrow M$ such that $\psi\varphi=1_M$.
A split mono is pure.
A pure map $M\rightarrow N$ is a split mono if $M/N$ is finitely presented
\cite[(5.2)]{HR2}.
A ring homomorphism $\varphi:R\rightarrow R'$ is said to be pure 
(resp.\ split) if it is pure (resp.\ a split mono) as an $R$-linear map.
A commutative ring $R$ of characteristic $p$ is said to be $F$-pure 
(resp.\ $F$-finite)
if the Frobenius map $F_R:R\rightarrow R$ given by $a\mapsto a^p$ is
pure (resp.\ finite) as a ring homomorphism.
$F$-purity was defined by Hochster--Roberts in 1970's \cite{HR} \cite{HR2}.
This notion is deeply connected with log-canonical 
singularity in characteristic zero \cite{Watanabe}.
An $\Bbb F_p$-scheme $X$
is said to be Frobenius split if $\O_{X^{(1)}}\rightarrow F_*\O_X$
splits as an $\O_{X^{(1)}}$-linear map \cite{MR}.
For the notation $X^{(1)}$, see \cite[(I.9.2)]{Jantzen}.
For an $F$-finite ring $R$ of characteristic $p$, $R$ is $F$-pure
if and only if $\Spec R$ is Frobenius split.
If a nonnegatively graded $F$-finite noetherian ring $R$ of characteristic
$p$ is $F$-pure, then $\Proj R$ is Frobenius split.

An $F$-finite noetherian ring $R$ of characteristic $p$ is said to be 
strongly $F$-regular if for any nonzerodivisor $a$ of $R$, the
$R^{(e)}$-linear map $aF^e: R^{(e)}\rightarrow R$ ($x\mapsto ax^{p^e}$) is
$R^{(e)}$-split \cite{HH}.
A strongly $F$-regular $F$-finite ring is $F$-regular in the sense of 
Hochster--Huneke \cite{HH2}, and hence it is Cohen--Macaulay normal
(\cite[(4.2)]{HH3} and \cite[(0.10)]{Velez}).
Under the same assumption as in Theorem~\ref{main.thm}, 
$S^G$ is strongly $F$-regular \cite{Hashimoto}, 
and in particular, Cohen--Macaulay.

\proof[Proof of Theorem~\ref{main.thm}] 
Let $\Gamma\rightarrow [G,G]$ be the universal covering.
Then $S$ has a good filtration as a $\Gamma$-module by
\cite[(3.1.3), (3.2.7)]{Donkin}, and $U$ is identified
with the unipotent radical of a Borel subgroup of $\Gamma$.
Thus without loss of generality, we may assume that $G$ is
semisimple and simply connected.

Let $\rho$ denote the half sum of positive roots.
Let $\St$ denote the first Steinberg module $\nabla_G((p-1)\rho)$.
For a $G$-module $W$ and $r\geq 0$, $W^{(r)}$ denotes the $r$th 
Frobenius twist of $W$ \cite[(I.9.10)]{Jantzen}.
Note that $W^{(r)}=W$ as an abelian group.
$w\in W$, viewed as an element of $W^{(r)}$ is sometimes
denoted by $w^{(r)}$.
If $R$ is a $G$-algebra, then $R^{(r)}$ is also a $G$-algebra, and
the $r$th Frobenius map $F^r:R^{(r)}\rightarrow R$ is a $G$-algebra map
\cite{Hashimoto}.

\begin{lemma}\label{s-split.thm}
There is a $(G,S^{(1)})$-linear splitting of the Frobenius map
$1\otimes F_S: \St\otimes S^{(1)}\rightarrow \St\otimes S$ given by 
$x\otimes s^{(1)}\mapsto x\otimes s^p$.
\end{lemma}

\proof The same proof as that of (4) in the proof of 
\cite[Theorem~6]{Hashimoto} works ($r=1$ and $d=0$ there).
\qed

Let $C:=k[G]^U$, where $U$ acts on $k[G]$ right regularly.
Then $C=\bigoplus_{\lambda\in X^+}\nabla_G(\lambda)$, and it 
is an $X^+$-graded $G$-algebra.

\begin{lemma}\label{c-split.thm}
There is a $(G,C^{(1)})$-linear splitting of the Frobenius map
$1\otimes F_C: \St\otimes C^{(1)}\rightarrow \St\otimes C$ given by
$x\otimes c^{(1)}\mapsto x\otimes c^p$.
\end{lemma}

\proof 
The product 
$\St\otimes \nabla_G(\lambda)^{(1)}\rightarrow
\nabla_G(p(\lambda+\rho)-\rho)$ is nonzero, since for
$x\in \St\setminus 0$ and 
$y\in \nabla_G(\lambda)\setminus 0$, $xy^{p}\neq 0$, 
as $C$ is an integral domain.
Since $\St\otimes\nabla_G(\lambda)^{(1)}\cong \nabla_G(p(\lambda+\rho)-\rho)$
\cite[(II.3.19)]{Jantzen} and $\End_G(\nabla_G(\mu))\cong k$ for each $\mu
\in X^+$, the product $\St\otimes\nabla_G(\lambda)^{(1)}\rightarrow
\nabla_G(p(\lambda+\rho)-\rho)$ is an isomorphism.
This shows that the product
$m:\St\otimes C^{(1)}\rightarrow C_{(p-1)\rho+pX^+}$ is an isomorphism,
where for a subset $\Lambda$ of 
$X^+$, $C_\Lambda:=\bigoplus_{\lambda\in\Lambda}\nabla_G(\lambda)$.

Thus
\[
\St\otimes C\xrightarrow{\pi}
\St\otimes C_{pX^+}
\xrightarrow{m'}C_{(p-1)\rho+pX^+}
\xrightarrow{m^{-1}}\St\otimes C^{(1)}
\]
is the splitting of $1\otimes F$, where $\pi$ is the projection, and
$m'$ is the product.
\qed

\begin{lemma}\label{sc-splitting.thm}
There exists some $(G,(S\otimes C)^{(1)})$-linear splitting
of the Frobenius map
$1\otimes F:\St\otimes (S\otimes C)^{(1)}\rightarrow \St\otimes S\otimes C$.
\end{lemma}

\proof Follows immediately from Lemma~\ref{s-split.thm} and 
Lemma~\ref{c-split.thm}.
\qed

\proof[Proof of Theorem~\ref{main.thm} (continued)]
Let $\Phi:\St\otimes S\otimes C\rightarrow \St\otimes (S\otimes C)^{(1)}$
be a $(G,(S\otimes C)^{(1)})$-linear splitting of $1\otimes F$ as in
Lemma~\ref{sc-splitting.thm}.
Set $A:=(S\otimes C)^G$, and consider the commutative
diagram of $(G,A^{(1)})$-modules
\[
\xymatrix{
\St\otimes A^{(1)} 
\ar@{^{(}->}[r]
\ar[d]^{1\otimes F} &
\St\otimes (S\otimes C)^{(1)}
\ar[r]^{\id}
\ar[d]^{1\otimes F} &
\St\otimes (S\otimes C)^{(1)}\\
\St\otimes A 
\ar@{^{(}->}[r] 
&
\St\otimes(S\otimes C)
\ar[ur]_{\Phi}.
}
\]
Applying the functor $\Hom_G(\St,?)$ to this, we get a commutative diagram
\[
\xymatrix{
A^{(1)} \ar[r]^{\id} \ar[d]^F &
A^{(1)} \ar[r]^{\id} \ar[d] &
A^{(1)} \\
A \ar[r] &
\Hom_G(\St,\St\otimes(S\otimes C)),
\ar[ur]
}
\]
of $A^{(1)}$-modules, see \cite[Proposition~1, {\bf 5}]{Hashimoto}.
Thus $F:A^{(1)}\rightarrow A$ splits, and $A=(S\otimes C)^G$ is $F$-pure.

Finally, as in the proof of \cite[(1.2)]{Grosshans}, 
$A=(S\otimes C)^G\cong S^U$ ($X=\Spec S$ and $H=U$ there).
Thus $S^U$ is $F$-pure.
\qed

\begin{corollary}
Under the same assumption as in Theorem~\ref{main.thm}, 
$\Proj S^U$ is Frobenius split.
\end{corollary}

\medskip
Acknowledgement: This research was done during the author's stay in
Tata Institute of Fundamental Research, Mumbai.
The author is grateful to TIFR and the host Professor V. Srinivas for
warm hospitality.
Special thanks are due to Professor V. Srinivas and Professor V. B. Mehta
for valuable discussion during the stay.

\end{document}